\documentstyle[twoside]{article}
\title{ On a generalized n-inner product \\
and the corresponding Cauchy - Schwarz inequality}

\author{Kostadin Ten\v{c}evski$^1$ and Risto Mal\v{c}eski$^2$}

\date{}
\begin{document}
\maketitle

\begin{abstract}

In this paper is defined an $n$-inner product of type 
$\langle {\bf a}_1,\cdots ,{\bf a}_n\vert {\bf b}_1\cdots {\bf b}_n\rangle $
where ${\bf a}_1,\cdots ,{\bf a}_n$, ${\bf b}_1,
\cdots ,{\bf b}_n$ are vectors from 
a vector space $V$. This definition generalizes the definition of Misiak
of $n$-inner product \cite{2}, 
such that in special case if we consider only such pairs of sets 
$\{ {\bf a}_1,\cdots ,{\bf a}_1\}$ and $\{ {\bf b}_1\cdots {\bf b}_n\}$ 
which differ for at most one vector, we obtain the definition of Misiak. 
The Cauchy-Schwarz inequality for this general type of $n$-inner 
product is proved and some applications are given. 
\medskip 


\end{abstract}

\section{Introduction}

A.Misiak \cite{2} has introduced an $n$-inner product by the following 
definition. 

\smallskip
\noindent {\bf Definition 1.1.} 
{\em Assume that $n$ is a positive integer and $V$ is a real
vector space such that $dim V\ge n$ and 
$(\bullet ,\bullet \vert \underbrace{ \bullet ,\cdots ,\bullet }_{n-1})$ 
be a real function defined on 
$\underbrace{V\times V \times \cdots \times V}_{n+1}$ such that }

i)\quad $({\bf x}_1,{\bf x}_1\vert {\bf x}_{2},\cdots ,{\bf x}_{n})\ge 0$, 
{\em for any 
${\bf x}_1,{\bf x}_2,\cdots ,{\bf x}_{n}\in V$ and 

\noindent $({\bf x}_1,{\bf x}_1\vert {\bf x}_{2},\cdots ,{\bf x}_{n})=0$ 
if and only if ${\bf x}_1,{\bf x}_2,\cdots ,{\bf x}_{n}$ 
are linearly dependent vectors; }

ii)\quad $({\bf a},{\bf b}\vert {\bf x}_{1},\cdots ,{\bf x}_{n-1})=
(\varphi ({\bf a}),\varphi ({\bf b})\vert \pi ({\bf x}_{1}),\cdots ,\pi 
({\bf x}_{n-1}))$, 

\noindent{\em for any ${\bf a},{\bf b},{\bf x}_1,\cdots ,{\bf x}_{n-1}\in V$ 
and for any bijections }
$$\pi :\{{\bf x}_1, \cdots ,{\bf x}_{n-1}\}\rightarrow \{{\bf x}_1, \cdots 
,{\bf x}_{n-1}\}\quad and \quad 
\varphi:\{{\bf a},{\bf b}\} \rightarrow \{{\bf a},{\bf b}\};$$

iii) {\em If $n>1$, then 
$({\bf x}_1,{\bf x}_1\vert {\bf x}_2, \cdots ,{\bf x}_{n})=
({\bf x}_2,{\bf x}_2\vert {\bf x}_1,{\bf x}_3, \cdots ,{\bf x}_{n})$, 

\noindent for any ${\bf x}_1,{\bf x}_2,\cdots ,{\bf x}_{n}\in V$; }

iv) \quad $(\alpha {\bf a},{\bf b}\vert {\bf x}_1, \cdots ,{\bf x}_{n-1})=
\alpha ({\bf a},{\bf b}\vert {\bf x}_1, \cdots ,{\bf x}_{n-1})$, 

\noindent {\em for any 
${\bf a}, {\bf b}, {\bf x}_1,\cdots ,{\bf x}_{n-1}\in V$ and any scalar 
$\alpha \in R$; }

v)\quad $({\bf a}+{\bf a}_1,{\bf b}\vert {\bf x}_1, \cdots ,{\bf x}_{n-1})=
({\bf a},{\bf b}\vert {\bf x}_1, \cdots ,{\bf x}_{n-1})+({\bf a}_1,{\bf 
b}\vert {\bf x}_1, \cdots ,{\bf x}_{n-1})$, 

\noindent {\em for any 
${\bf a},{\bf b},{\bf a}_1,{\bf x}_1,\cdots ,{\bf x}_{n-1}\in V$. 

Then 
$(\bullet ,\bullet \vert \underbrace{ \bullet ,\cdots ,\bullet }_{n-1})$ 
is called $n$-inner product and 
$(V,(\bullet ,\bullet \vert \underbrace{ \bullet ,\cdots ,\bullet }_{n-1}))$
is called $n$-prehilbert space.} 

If $n=1$, then Definition 1.1 reduces to the ordinary inner product. 

This $n$-inner product induces an $n$-norm (\cite{2}) by 
$$\Vert {\bf x}_1,\cdots ,{\bf x}_n\Vert = 
\sqrt{ ({\bf x}_1,{\bf x}_1\vert {\bf x}_2, \cdots ,{\bf x}_n)}.$$
In the next section 
we introduce a more general and more convenient 
definition of $n$-inner product 
and prove the corresponding Cauchy-Schwarz inequality. 
In the last section some related results are given. 

Although in this paper are considered only real vector spaces, 
the results of this paper can easily be generalized for 
the complex vector spaces.

\section{n-inner product and the Cauchy-Schwarz inequality}

First we give the following definition of $n$-inner product. 

\smallskip
\noindent {\bf Definition 2.1.} {\em Assume that $n$ is a positive integer, 
$V$ is a real vector space such that 
$\dim V\ge n$ and 
$\langle \bullet,\cdots ,\bullet\vert \bullet,\cdots ,\bullet \rangle$ be a 
real function
on $V^{2n}$ such that }
\par 
i) \quad $$\langle {\bf a}_{1},\cdots ,{\bf a}_{n}\vert 
{\bf a}_{1},\cdots ,{\bf a}_{n}\rangle > 0\leqno{(2.1)}$$ 
{\em if ${\bf a}_{1},\cdots ,{\bf a}_{n}$ are linearly independent vectors,} 

ii) \quad 
$$\langle {\bf a}_{1},\cdots ,{\bf a}_{n}\vert {\bf b}_{1},\cdots ,{\bf 
b}_{n}\rangle = 
\langle {\bf b}_{1},\cdots ,{\bf b}_{n}\vert {\bf a}_{1},\cdots ,{\bf 
a}_{n}\rangle
\leqno{(2.2)}$$ 
{\em for any }
${\bf a}_{1},\cdots ,{\bf a}_{n},{\bf b}_{1},\cdots ,{\bf b}_{n}\in V$, 

iii) \quad 
$$\langle \lambda {\bf a}_{1},\cdots ,{\bf a}_{n}\vert 
{\bf b}_{1},\cdots ,{\bf b}_{n}\rangle = \lambda 
\langle {\bf a}_{1},\cdots ,{\bf a}_{n}\vert {\bf b}_{1},\cdots ,{\bf 
b}_{n}\rangle
\leqno{(2.3)}$$ 
{\em for any scalar $\lambda \in R$ and any }
${\bf a}_{1},\cdots ,{\bf a}_{n},{\bf b}_{1},\cdots ,{\bf b}_{n}\in V$, 

iv) \quad 
$$\langle {\bf a}_{1},\cdots ,{\bf a}_{n}\vert {\bf b}_{1},\cdots ,{\bf 
b}_{n}\rangle = -
\langle {\bf a}_{\sigma (1)},\cdots ,{\bf a}_{\sigma (n)}
\vert {\bf b}_{1},\cdots ,{\bf b}_{n}\rangle\leqno{(2.4)}$$ 
{\em for any odd permutation 
$\sigma $ in the set $\{ 1,\cdots ,n\}$ and any} 
${\bf a}_{1},\cdots ,{\bf a}_{n},{\bf b}_{1},\cdots ,{\bf b}_{n}\in V$, 

v) \quad 
$$\langle {\bf a}_{1}+{\bf c},{\bf a}_{2},\cdots ,{\bf a}_{n}
\vert {\bf b}_{1},\cdots ,{\bf b}_{n}\rangle = 
\langle {\bf a}_{1},{\bf a}_{2},\cdots ,{\bf a}_{n}
\vert {\bf b}_{1},\cdots ,{\bf b}_{n}\rangle + 
\langle {\bf c},{\bf a}_{2},\cdots ,{\bf a}_{n}
\vert {\bf b}_{1},\cdots ,{\bf b}_{n}\rangle \leqno{(2.5)}$$
{\em for any} 
${\bf a}_{1},\cdots ,{\bf a}_{n},{\bf b}_{1},\cdots ,
{\bf b}_{n},{\bf c}\in V$, 
\par 
vi) \quad {\em if } 
$$\langle {\bf a}_{1},{\bf b}_{1},\cdots ,{\bf b}_{i-1},{\bf b}_{i+1},
\cdots ,{\bf b}_{n}\vert {\bf b}_{1},\cdots ,{\bf b}_{n}\rangle = 0
\leqno{(2.6')}$$
{\em for each $i\in \{1,2,\cdots ,n\}$, then }
$$\langle {\bf a}_{1},\cdots ,{\bf a}_{n}\vert {\bf b}_{1},\cdots ,{\bf 
b}_{n}\rangle =0\leqno{(2.6'')}$$ 
{\em for arbitrary vectors} ${\bf a}_{2},\cdots ,{\bf a}_{n}$. 

{\em Then the function 
$\langle \bullet,\cdots ,\bullet\vert \bullet,\cdots ,\bullet \rangle$ 
is called 
$n$-inner product and the pair 
$(V,\langle \bullet,\cdots ,\bullet\vert \bullet,\cdots ,\bullet \rangle)$ 
is called $n$-prehilbert space.} 

We give some consequences from the conditions i)-vi) of Definition 2.1. 

{From} (2.4) it follows that if two of the vectors 
${\bf a}_{1},\cdots ,{\bf a}_{n}$ are equal, then 
$\langle {\bf a}_{1},\cdots ,{\bf a}_{n}\vert {\bf b}_{1},\cdots ,{\bf 
b}_{n}\rangle=0.$

{From} (2.3) it follows that 
$$\langle {\bf a}_{1},\cdots ,{\bf a}_{n}\vert {\bf b}_{1},\cdots ,{\bf 
b}_{n}\rangle=0$$
if there exists $i$ such that ${\bf a}_{i}=0$. 

{From} (2.4) and (2.2) it follows more generally that 

iv$'$)\quad 
$$\langle {\bf a}_{1},\cdots ,{\bf a}_{n}\vert {\bf b}_{1},\cdots ,{\bf 
b}_{n}\rangle = 
(-1)^{sgn (\pi )+ sgn (\tau )} 
\langle {\bf a}_{\pi (1)},\cdots ,{\bf a}_{\pi (n)}\vert 
{\bf b}_{\tau (1)},\cdots ,{\bf b}_{\tau (n)}\rangle $$
{\em for any permutations 
$\pi $ and $\tau $ on $\{1,\cdots ,n\}$ 
and } 
${\bf a}_{1},\cdots ,{\bf a}_{n},{\bf b}_{1},\cdots ,{\bf b}_{n}\in V$.

{From} (2.3), (2.4) and (2.5) it follows that 
$$\langle {\bf a}_{1},\cdots ,{\bf a}_{n}\vert {\bf b}_{1},\cdots ,{\bf 
b}_{n}\rangle=0$$
if ${\bf a}_{1},\cdots ,{\bf a}_{n}$ are linearly dependent vectors. Thus 
i) can be replaced by 

i$'$) \quad 
$\langle {\bf a}_{1},\cdots ,{\bf a}_{n}\vert {\bf a}_{1},\cdots ,{\bf 
a}_{n}\rangle\ge 0$

\noindent {\em for any 
${\bf a}_{1},\cdots ,{\bf a}_{n}\in V$ and 
$\langle {\bf a}_{1},\cdots ,{\bf a}_{n}\vert {\bf a}_{1},\cdots ,{\bf 
a}_{n}\rangle=0$
if and only if 
${\bf a}_{1},\cdots ,{\bf a}_{n}$ are linearly dependent vectors.} 

Note that the $n$-inner product on $V$ induces an $n$-normed space by 
$$\Vert {\bf x}_{1},\cdots ,{\bf x}_{n}\Vert = \sqrt{
\langle {\bf  x}_{1},\cdots  ,{\bf  x}_{n}\vert  {\bf  x}_{1},\cdots ,{\bf 
x}_{n}\rangle} ,$$
and it is the same norm induced by the Definition 1.1. 

In special case if we consider only such pairs of sets 
${\bf a}_1,\cdots ,{\bf a}_1$ and ${\bf b}_1\cdots {\bf b}_n$ which differ 
for at most one vector, for example ${\bf a}_1={\bf a}$, ${\bf b}_1={\bf b}$
and ${\bf a}_2={\bf b}_2={\bf x}_1,
\cdots ,{\bf a}_n={\bf b}_n={\bf x}_{n-1}$, then by putting 
$$({\bf a},{\bf b}\vert {\bf x}_1,\cdots ,{\bf x}_{n-1})=
\langle {\bf a},{\bf x}_1,\cdots ,{\bf x}_{n-1}\vert 
{\bf b},{\bf x}_1,\cdots ,{\bf x}_{n-1}\rangle$$
we obtain an $n$-inner product according to Definition 1.1 of Misiak. 
Indeed, the 
conditions i), iv) and v) are triavially satisfied. The condition ii) 
is satisfied for arbitrary permutation $\pi $, because according 
to iv$'$) 
$$\langle {\bf a}_{1},\cdots ,{\bf a}_{n}\vert {\bf b}_{1},\cdots ,{\bf 
b}_{n}\rangle = 
\langle {\bf a}_1,{\bf a}_{\pi (2)},\cdots ,{\bf a}_{\pi (n)}\vert 
{\bf b}_1,{\bf b}_{\pi (2)},\cdots ,{\bf b}_{\pi (n)}\rangle $$
for any permutation 
$\pi :\{ 2,3,\cdots, n\}\rightarrow \{ 2,3,\cdots, n\}$. Similarly 
the condition iii) is satisfied. Moreover, in this special case of the 
Definition 2.1 we do not have any restriction of Definition 1.1. 
For example, then the condition vi) does not say anything. Namely, 
if ${\bf a}_1\notin \{ {\bf b}_1,\cdots ,{\bf b}_n\}$, then the vectors 
${\bf a}_2,\cdots ,{\bf a}_n$ must be from the set 
$\{ {\bf b}_1,\cdots ,{\bf b}_n\}$, and $(2.6'')$ is satisfied because 
the assumption $(2.6')$ is satisfied. If 
${\bf a}_1\in \{ {\bf b}_1,\cdots ,{\bf b}_n\}$, for example 
${\bf a}_1={\bf b}_j$, then $(2.6')$ implies that 
$\langle {\bf b}_j,{\bf b}_1,\cdots ,{\bf b}_{j-1},{\bf b}_{j+1},\cdots ,
{\bf b}_n \vert {\bf b}_1,\cdots ,{\bf b}_n\rangle =0$, and it i possible 
only if ${\bf b}_1,\cdots ,{\bf b}_n$ are linearly dependent vectors. 
But, then $(2.6'')$ is satisfied. 
Thus the Definition 2.1 generalizes Definition 1.1.

Now we give the following example of the $n$-inner product. 

\smallskip
\noindent {\bf Example 2.1.} We refer to the classical known example, as an 
$n$-inner product according to Definition 2.1. 
Let $V$ be a space with inner product $\langle \vert \rangle$. 
Then 
$$\langle {\bf a}_{1},\cdots ,{\bf a}_{n}\vert {\bf b}_{1},\cdots ,{\bf 
b}_{n}\rangle = 
\left \vert \matrix{ 
\langle {\bf a}_{1}\vert {\bf b}_{1}\rangle&\langle {\bf a}_{1}\vert {\bf 
b}_{2}\rangle&\cdots &
\langle {\bf a}_{1}\vert {\bf b}_{n}\rangle\cr 
\langle {\bf a}_{2}\vert {\bf b}_{1}\rangle&\langle {\bf a}_{2}\vert {\bf 
b}_{2}\rangle&\cdots &
\langle {\bf a}_{2}\vert {\bf b}_{n}\rangle\cr 
\cdot & & & \cr
\cdot & & & \cr
\cdot & & & \cr
\langle {\bf a}_{n}\vert {\bf b}_{1}\rangle&\langle {\bf a}_{n}\vert {\bf 
b}_{2}\rangle&\cdots &
\langle {\bf a}_{n}\vert {\bf b}_{n}\rangle\cr }\right \vert $$
satisfies the conditions i) - vi) and hence it defines an 
$n$-inner product on 
$V$. The conditions i) - v) are trivial, and we will prove vi). If 
${\bf b}_{1},\cdots ,{\bf b}_{n}$ are linearly independent vectors and 
$$\langle {\bf a}_{1},{\bf b}_{1},\cdots ,{\bf b}_{i-1},{\bf b}_{i+1},
\cdots ,{\bf b}_{n}\vert {\bf b}_{1},\cdots ,{\bf b}_{n}\rangle \equiv 
(-1)^{i-1}\langle {\bf b}_{1},\cdots ,{\bf b}_{i-1},{\bf a}_1,{\bf b}_{i+1},
\cdots ,{\bf b}_{n}\vert {\bf b}_{1},\cdots ,{\bf b}_{n}\rangle \equiv $$
$$\equiv (-1)^{i-1}
\left \vert \matrix{ 
\langle {\bf b}_{1}\vert {\bf b}_{1}\rangle&\langle {\bf b}_{1}\vert {\bf 
b}_{2}\rangle&\cdots & \langle {\bf b}_{1}\vert {\bf b}_{n}\rangle\cr 
\langle {\bf b}_{2}\vert {\bf b}_{1}\rangle&\langle {\bf b}_{2}\vert {\bf 
b}_{2}\rangle&\cdots&\langle {\bf b}_{2}\vert {\bf b}_{n}\rangle\cr 
\cdots & \cdots & \cdots & \cdots \cr
\langle {\bf a}_{1}\vert {\bf b}_{1}\rangle&\langle {\bf a}_{1}\vert {\bf 
b}_{2}\rangle&\cdots &\langle {\bf a}_{1}\vert {\bf b}_{n}\rangle\cr 
\cdots & \cdots & \cdots & \cdots \cr
\langle {\bf b}_{n}\vert {\bf b}_{1}\rangle&\langle {\bf b}_{n}\vert {\bf 
b}_{2}\rangle&\cdots &\langle {\bf b}_{n}\vert {\bf b}_{n}\rangle\cr }
\right \vert =0,$$
then the vector 
$$(\langle {\bf a}_{1}\vert {\bf b}_{1}\rangle ,\langle {\bf a}_{1}\vert {\bf 
b}_{2}\rangle ,\cdots ,\langle {\bf a}_{1}\vert {\bf b}_{n}\rangle )\in R^n$$
is a linear combination of 
$$(\langle {\bf b}_{1}\vert {\bf b}_{1}\rangle ,
\cdots ,\langle {\bf b}_{1}\vert {\bf b}_{n}\rangle ), \cdots ,
(\langle {\bf b}_{i-1}\vert {\bf b}_{1}\rangle ,
\cdots ,\langle {\bf b}_{i-1}\vert {\bf b}_{n}\rangle ), $$
$$(\langle {\bf b}_{i+1}\vert {\bf b}_{1}\rangle ,
\cdots ,\langle {\bf b}_{i+1}\vert {\bf b}_{n}\rangle ), \cdots ,
(\langle {\bf b}_{n}\vert {\bf b}_{1}\rangle ,
\cdots ,\langle {\bf b}_{n}\vert {\bf b}_{n}\rangle ). $$
Since this is true for each $i\in \{ 1,2,\cdots ,n\}$, it must be 
$\langle {\bf a}_{1}\vert {\bf b}_{1}\rangle =\cdots =
\langle {\bf a}_{1}\vert {\bf b}_{n}\rangle =0$. Hence 
$$\langle {\bf a}_{1},\cdots ,{\bf a}_{n}\vert {\bf b}_{1},\cdots ,{\bf 
b}_{n}\rangle =0 $$ 
for arbitrary ${\bf a}_{2,},\cdots ,{\bf a}_{n}$. 

Note that the inner product defined by 
$$\langle {\bf a}_{1}\wedge \cdots \wedge {\bf a}_{n}\vert 
{\bf b}_{1}\wedge \cdots \wedge {\bf b}_{n}\rangle = 
\left \vert \matrix{ 
\langle {\bf a}_{1}\vert {\bf b}_{1}\rangle&\langle {\bf a}_{1}\vert {\bf 
b}_{2}\rangle&\cdots &
\langle {\bf a}_{1}\vert {\bf b}_{n}\rangle\cr 
\langle {\bf a}_{2}\vert {\bf b}_{1}\rangle&\langle {\bf a}_{2}\vert {\bf 
b}_{2}\rangle&\cdots &
\langle {\bf a}_{2}\vert {\bf b}_{n}\rangle\cr 
\cdot & & & \cr
\cdot & & & \cr
\cdot & & & \cr
\langle {\bf a}_{n}\vert {\bf b}_{1}\rangle&\langle {\bf a}_{n}\vert {\bf 
b}_{2}\rangle&\cdots &
\langle {\bf a}_{n}\vert {\bf b}_{n}\rangle\cr }\right \vert $$
can uniquely be prolonged to ordinary inner product over the space 
$\Lambda _{n}(V)$ of $n$-forms over V \cite{1}. Indeed, if 
$\{ {\bf e}_{i}\}_{i\in I}$, $I$ an index set, 
is an orthonormal basis of $(V,\langle *\vert *\rangle)$, then 
$$ \langle {\bf e}_{i_{1}}\wedge \cdots \wedge {\bf e}_{i_{n}}\vert 
{\bf e}_{j_{1}}\wedge \cdots \wedge {\bf e}_{j_{n}}\rangle = 
\delta ^{i_{1}\cdots i_{n}}_{j_{1}\cdots j_{n}} $$
where the expression 
$\delta ^{i_{1}\cdots i_{n}}_{j_{1}\cdots j_{n}}$ is equal to 1 or -1 
if $\{ i_1,\cdots ,i_n\} = \{ j_1,\cdots ,j_n\}$ with different 
$i_1,\cdots ,i_n$ and additionally the permutation 
${i_1\; i_2\; \cdots \;i_n\choose j_1\; j_2\; \cdots \;j_n}$ 
is even or odd respectively, and where the above expression is 0 otherwise. 
It implies an inner 
product over $\Lambda _{n}(V)$. 

Before we prove the next theorem, we give the following remarks assuming 
that $\dim V>n$. 

Let ${\bf b}_{1},\cdots ,{\bf b}_{n}$ are linearly independent vectors. If 
a vector ${\bf a}$ is such that 
$$\langle {\bf a},{\bf b}_{1},\cdots ,{\bf b}_{i-1},{\bf b}_{i+1},
\cdots ,{\bf b}_{n}\vert {\bf b}_{1},\cdots ,{\bf b}_{n}\rangle =0,
\quad (1\le i\le n)$$
then we say that the vector ${\bf a}$ is orthogonal to the subspace 
generated by ${\bf b}_{1},\cdots ,{\bf b}_{n}$. Note that the set of 
orthogonal vectors to this $n$-dimensional subspace is a vector subspace of 
$V$, and the orthogonality of ${\bf a}$ to the considered vector subspace 
is invariant of the base vectors ${\bf b}_{1},\cdots ,{\bf b}_{n}$. 
If ${\bf x}$ is an arbitrary vector, then there exist unique 
$\lambda _1,\cdots ,\lambda _n\in R$ such that 
${\bf x}-\lambda_1 {\bf b}_1-\cdots -\lambda_n {\bf b}_n$ is orthogonal to 
the vector subspace generated by ${\bf b}_{1},\cdots ,{\bf b}_{n}$. Namely, 
the orthogonality conditions 
$$\langle {\bf b}_{1},\cdots ,{\bf b}_{i-1},
{\bf x}-\lambda_1 {\bf b}_1-\cdots -\lambda_n {\bf b}_n,
{\bf b}_{i+1},\cdots ,{\bf b}_{n}\vert 
{\bf b}_{1},\cdots ,{\bf b}_{n}\rangle =0,
\quad (1\le i\le n)$$
have unique solutions 
$$\lambda _i={
\langle {\bf b}_{1},\cdots ,{\bf b}_{i-1},{\bf x}, {\bf b}_{i+1},
\cdots ,{\bf b}_{n}\vert {\bf b}_{1},\cdots ,{\bf b}_{n}\rangle \over 
\langle {\bf b}_{1},\cdots ,{\bf b}_{n}\vert {\bf b}_{1},\cdots ,{\bf b}_{n}
\rangle }, \quad (1\le i\le n).$$
Hence each vector ${\bf x}$ can uniquely be decomposed as 
${\bf x}=\lambda _1{\bf b}_1+\cdots +\lambda _n{\bf b}_n+{\bf c}$, 
where the vector ${\bf c}$ is orthogonal to the vector subspace generated 
by ${\bf b}_{1},\cdots ,{\bf b}_{n}$. According to this definition, 
the condition vi) of Definition 2.1 says that if the vector 
${\bf a}_1$ is orthogonal to vector subspace generated by 
${\bf b}_{1},\cdots ,{\bf b}_{n}$, then $(2.6'')$ holds for arbitrary 
vectors ${\bf a}_2,\cdots ,{\bf a}_n$. 

Now we prove the Cauchy-Schwarz inequality as 
a consequence of Definition 2.1. 

\smallskip
\noindent {\bf Theorem 2.1.} {\em If} 
$\langle \bullet, \cdots ,\bullet \vert \bullet, 
\cdots ,\bullet \rangle$ {\em is an 
$n$-inner product on $V$, then the following inequality} 
$$
\langle {\bf a}_{1},\cdots ,{\bf a}_{n}\vert {\bf b}_{1},\cdots ,{\bf 
b}_{n}\rangle^{2} \le 
\langle {\bf a}_{1},\cdots ,{\bf a}_{n}\vert {\bf a}_{1},\cdots ,{\bf 
a}_{n}\rangle
\langle {\bf b}_{1},\cdots ,{\bf b}_{n}\vert {\bf b}_{1},\cdots ,{\bf 
b}_{n}\rangle ,
\leqno{(2.7)}$$
{\em is true for any vectors 
${\bf a}_{1},\cdots ,{\bf a}_{n},{\bf b}_{1},\cdots ,{\bf b}_{n}\in V$. 
Moreover, equality holds if and only if at least one of the following 
conditions is satisfied }

i) \quad {\em the vectors ${\bf a}_{1},{\bf a}_{2},\cdots ,{\bf a}_{n}$ are 
linearly dependent,}
\par 
ii) \quad {\em the vectors ${\bf b}_{1},{\bf b}_{2},\cdots ,{\bf b}_{n}$ are 
linearly dependent,}
\par 
iii) \quad {\em the vectors ${\bf a}_{1},{\bf a}_{2},\cdots ,{\bf a}_{n}$ 
and ${\bf b}_{1},{\bf b}_{2},\cdots ,{\bf b}_{n}$ generate 
the same vector subspace of dimension $n$. }

\smallskip
\noindent 
{\em Proof.} If ${\bf a}_{1},\cdots ,{\bf a}_{n}$ are linearly dependent 
vectors or ${\bf b}_{1},\cdots ,{\bf b}_{n}$ are linearly dependent 
vectors, then both sides of (2.7) are zero and hence equality holds. Thus, 
suppose that ${\bf a}_{1},\cdots ,{\bf a}_{n}$ and also 
${\bf b}_{1},\cdots ,{\bf b}_{n}$ are linearly independent vectors. 
Note that the inequality (2.7) does not depend on the choice 
of the basis ${\bf a}_{1},\cdots ,{\bf a}_{n}$ of the subspace
generated by these $n$ vectors. 
Indeed, each vector row operation preserves the inequality 
(2.7), 
because both sides are invariant or both sides are multiplied by 
a positive real scalar after any elementary vector row operation. 
We assume that $\dim V>n$, because if $\dim V=n$, then 
the theorem is obviously satisfied. 

Let $\Sigma $ be a space generated by the vectors 
${\bf a}_{1},\cdots ,{\bf a}_{n}$ and $\Sigma ^{*}$ be the orthogonal 
subspace to $\Sigma $. Let us decompose the vectors 
${\bf b}_{i}$ as ${\bf b}_{i}={\bf c}_{i}+{\bf d}_{i}$ where 
${\bf c}_{i}\in \Sigma $ and ${\bf d}_{i}\in \Sigma ^{*}$. 
Thus 
$$ {\bf b}_{i} = \sum _{j=1}^{n}P_{ij}{\bf a}_{j} + {\bf d}_{i}, \qquad 
(1\le i\le n)$$
$$\langle {\bf a}_{1},\cdots ,{\bf a}_{n}\vert {\bf b}_{1},\cdots ,{\bf 
b}_{n}\rangle
= \langle {\bf a}_{1},\cdots ,{\bf a}_{n}\vert 
\sum _{j_{1}=1}^{n}P_{1j_{1}}{\bf a}_{j_{1}} + {\bf d}_{1}, \cdots ,
\sum _{j_{n}=1}^{n}P_{nj_{n}}{\bf a}_{j_{n}} + {\bf d}_{n}\rangle$$
$$ = \sum _{j_{1}=1}^{n} \cdots \sum _{j_{n}=1}^{n} 
P_{1j_{1}}P_{2j_{2}}\cdots P_{nj_{n}}
\langle {\bf a}_{1},\cdots ,{\bf a}_{n}\vert 
{\bf a}_{j_{1}},\cdots ,{\bf a}_{j_{n}}\rangle$$
$$ = \sum _{j_{1}=1}^{n} \cdots \sum _{j_{n}=1}^{n} 
P_{1j_{1}}P_{2j_{2}}\cdots P_{nj_{n}} (-1)^{sgn \sigma }
\langle {\bf a}_{1},\cdots ,{\bf a}_{n}\vert 
{\bf a}_{1},\cdots ,{\bf a}_{n}\rangle$$
$$ = det P \cdot 
\langle {\bf a}_{1},\cdots ,{\bf a}_{n}\vert 
{\bf a}_{1},\cdots ,{\bf a}_{n}\rangle$$
where we used the conditions ii) - vi) from definition 2.1 
and we denoted by $P$ 
the matrix with entries $P_{ij}$, and 
$\sigma ={1\; 2\; \cdots \;n\choose j_1 j_2 \cdots j_n}$. 

If $det P=0$, then the left side of (2.7) is 0, the right side is positive 
and hence the inequality (2.7) is true. So, let us suppose that $P$ is 
a non-singular matrix and $Q=P^{-1}$. Now the inequality (2.7) is 
equivalent to 
$$(det P)^{2} 
\langle {\bf a}_{1},\cdots ,{\bf a}_{n}\vert {\bf a}_{1},\cdots ,{\bf 
a}_{n}\rangle^{2} \le 
\langle {\bf a}_{1},\cdots ,{\bf a}_{n}\vert {\bf a}_{1},\cdots ,{\bf 
a}_{n}\rangle
\langle {\bf b}_{1},\cdots ,{\bf b}_{n}\vert {\bf b}_{1},\cdots ,{\bf 
b}_{n}\rangle, $$
$$\langle {\bf a}_{1},\cdots ,{\bf a}_{n}\vert {\bf a}_{1},\cdots ,{\bf 
a}_{n}\rangle\le 
\langle {\bf b}'_{1},\cdots ,{\bf b}'_{n}\vert {\bf b}'_{1},\cdots ,{\bf 
b}'_{n}\rangle ,\leqno{(2.8)}$$
where ${\bf b}'_{i}=\sum\limits _{j=1}^{n}Q_{ij}{\bf b}_{j}$,\quad 
$(1\le i\le n)$. Note that ${\bf b}'_{i}$ decomposes as 
$$ {\bf b}'_{i} = \sum _{j=1}^{n}Q_{ij}\Bigl (
\sum _{l=1}^{n}P_{jl}{\bf a}_{l}+
{\bf d}_{j}\Bigr ) = {\bf a}_{i} + {\bf d}'_{i}$$
where 
${\bf d}'_{i}=\sum\limits _{j=1}^{n}Q_{ij}{\bf d}_{j}\in {\Sigma }^{*}$.
Now we will prove (2.8), i.e. 
$$\langle {\bf a}_{1},\cdots ,{\bf a}_{n}\vert {\bf a}_{1},\cdots ,{\bf 
a}_{n}\rangle\le 
\langle {\bf a}_{1}+{\bf d}'_{1},\cdots ,{\bf a}_{n}+{\bf d}'_{n}\vert 
{\bf a}_{1}+{\bf d}'_{1},\cdots ,
{\bf a}_{n}+{\bf d}'_{n}\rangle\leqno{(2.9)}$$
and equality holds if and only if ${\bf b}'_{1}={\bf a}_{1}$, $\cdots $, 
${\bf b}'_{n}={\bf a}_{n}$, i.e., 
${\bf d}'_{1}=\cdots ={\bf d}'_{n}=0$. More precisely, we will prove that 
(2.9) is true for at least one basis ${\bf a}_{1},\cdots ,{\bf a}_{n}$
of $\Sigma$. 

Using (2.5) and (2.2) we obtain 
$$\langle {\bf a}_{1}+{\bf d}'_{1},\cdots ,{\bf a}_{n}+{\bf d}'_{n}\vert 
{\bf a}_{1}+{\bf d}'_{1},\cdots ,{\bf a}_{n}+{\bf d}'_{n}\rangle $$
$$= \langle {\bf a}_{1},{\bf a}_{2}+{\bf d}'_{2}, 
\cdots ,{\bf a}_{n}+{\bf d}'_{n}\vert 
{\bf a}_{1},{\bf a}_{2}+{\bf d}'_{2},
\cdots ,{\bf a}_{n}+{\bf d}'_{n}\rangle $$
$$+ \langle {\bf d}'_{1},{\bf a}_{2}+{\bf d}'_{2}, 
\cdots ,{\bf a}_{n}+{\bf d}'_{n}\vert 
{\bf d}'_{1},{\bf a}_{2}+{\bf d}'_{2},
\cdots ,{\bf a}_{n}+{\bf d}'_{n}\rangle $$
$$+2\langle {\bf a}_{1},{\bf a}_{2}+{\bf d}'_{2},\cdots ,{\bf a}_{n}+
{\bf d}'_{n}\vert {\bf d}'_{1},{\bf a}_{2}+{\bf d}'_{2},
\cdots ,{\bf a}_{n}+{\bf d}'_{n}\rangle $$
$$ = \langle {\bf a}_{1},{\bf a}_{2},{\bf a}_{3}+{\bf d}'_{3}
,\cdots ,{\bf a}_{n}+{\bf d}'_{n}\vert 
{\bf a}_{1},{\bf a}_{2},{\bf a}_{3}+{\bf d}'_{3}
,\cdots ,{\bf a}_{n}+{\bf d}'_{n}\rangle$$
$$+
\langle {\bf a}_{1},{\bf d}'_{2},{\bf a}_{3}+{\bf d}'_{3}
,\cdots ,{\bf a}_{n}+{\bf d}'_{n}\vert 
{\bf a}_{1},{\bf d}'_{2},{\bf a}_{3}+{\bf d}'_{3}
,\cdots ,{\bf a}_{n}+{\bf d}'_{n}\rangle$$
$$+
\langle {\bf d}'_{1},{\bf a}_{2}+{\bf d}'_{2}
,\cdots ,{\bf a}_{n}+{\bf d}'_{n}\vert 
{\bf d}'_{1},{\bf a}_{2}+{\bf d}'_{2}
,\cdots ,{\bf a}_{n}+{\bf d}'_{n}\rangle $$
$$+2\langle {\bf a}_{1},{\bf a}_{2}+{\bf d}'_{2},\cdots ,{\bf a}_{n}+
{\bf d}'_{n}\vert {\bf d}'_{1},{\bf a}_{2}+{\bf d}'_{2},
\cdots ,{\bf a}_{n}+{\bf d}'_{n}\rangle $$
$$+2\langle {\bf a}_{1},{\bf a}_{2},{\bf a}_{3}+{\bf d}'_{3},\cdots ,
{\bf a}_{n}+{\bf d}'_{n}\vert {\bf a}_{1},
{\bf d}'_{2},{\bf a}_{3}+{\bf d}'_{3},
\cdots ,{\bf a}_{n}+{\bf d}'_{n}\rangle $$
$$ = \cdots = 
\langle {\bf a}_{1},\cdots ,{\bf a}_{n}\vert 
{\bf a}_{1},\cdots ,{\bf a}_{n}\rangle $$
$$+\langle {\bf a}_{1},\cdots ,{\bf a}_{n-1},{\bf d}'_{n}\vert 
{\bf a}_{1},\cdots ,{\bf a}_{n-1},{\bf d}'_{n}\rangle $$
$$+\cdots + 
\langle {\bf a}_{1},{\bf d}'_{2},\cdots ,{\bf a}_{n}+{\bf d}'_{n}\vert 
{\bf a}_{1},{\bf d}'_{2},\cdots ,{\bf a}_{n}+{\bf d}'_{n}\rangle$$
$$+
\langle {\bf d}'_{1},{\bf a}_{2}+{\bf d}'_{2}
,\cdots ,{\bf a}_{n}+{\bf d}'_{n}\vert 
{\bf d}'_{1},{\bf a}_{2}+{\bf d}'_{2}
,\cdots ,{\bf a}_{n}+{\bf d}'_{n}\rangle + S, $$
where 
$$S=
2\langle {\bf a}_{1},{\bf a}_{2}+{\bf d}'_{2},\cdots ,{\bf a}_{n}+
{\bf d}'_{n}\vert {\bf d}'_{1},{\bf a}_{2}+{\bf d}'_{2},
\cdots ,{\bf a}_{n}+{\bf d}'_{n}\rangle $$
$$+2\langle {\bf a}_{1},{\bf a}_{2},{\bf a}_{3}+{\bf d}'_{3},\cdots ,
{\bf a}_{n}+{\bf d}'_{n}\vert {\bf a}_{1},
{\bf d}'_{2},{\bf a}_{3}+{\bf d}'_{3},
\cdots ,{\bf a}_{n}+{\bf d}'_{n}\rangle $$
$$+\cdots + 
2\langle {\bf a}_{1},{\bf a}_{2},\cdots ,{\bf a}_{n-1},
{\bf a}_{n}\vert {\bf a}_{1},
{\bf a}_{2},\cdots ,{\bf a}_{n-1},{\bf d}'_{n}\rangle $$
We can change the basis ${\bf a}_{1},\cdots ,{\bf a}_{n}$ of $\Sigma$ 
such that the sum $S$ vanish. Indeed, if we replace 
${\bf a}_{1}$ by $\lambda {\bf a}_{1}$ we can choose almost always 
a scalar $\lambda $ such that $S=0$. The other cases can be considered 
by another analogous linear transformations. Thus without loss of 
generality we can put $S=0$. 

According to i$'$) the inequality (2.9) 
is true and equality holds if and only 
if the following sets of vectors 
$\{{\bf a}_{1},\cdots ,{\bf a}_{n-1},{\bf d}'_{n}\}$, $\cdots$, 
$\{{\bf a}_{1},{\bf d}'_{2},{\bf a}_{3}+{\bf d}'_{3},\cdots ,
{\bf a}_{n}+{\bf d}'_{n}\}$, 
$\{{\bf d}'_{1},{\bf a}_{2}+{\bf d}'_{2}
,\cdots ,{\bf a}_{n}+{\bf d}'_{n}\}$ are linearly dependent. This is 
satisfied if and only if 
${\bf d}'_{1}={\bf d}'_{2}=\cdots ={\bf d}'_{n}=0$, i.e. if and only if 
${\bf b}'_{1}={\bf a}_{1}$, $\cdots $, ${\bf b}'_{n}={\bf a}_{n}$. 
{\hspace* {\fill} { $\Box $}}

\section{Some applications}

Let $\Sigma _{1}$ and $\Sigma _{2}$ are two subspaces of $V$ 
of dimension $n$. We define 
angle $\varphi $ between $\Sigma _{1}$ and $\Sigma _{2}$ by 
$$ \cos \varphi = {
\langle {\bf a}_1,\cdots ,{\bf a}_n\vert {\bf b}_1,\cdots ,{\bf b}_n\rangle 
\over 
\Vert {\bf a}_1,\cdots ,{\bf a}_n\Vert \cdot 
\Vert {\bf b}_1,\cdots ,{\bf b}_n\Vert },\leqno{(3.1)}$$
where ${\bf a}_1,\cdots ,{\bf a}_n$ are linearly independent vectors of 
$\Sigma_1$, ${\bf b}_1,\cdots ,{\bf b}_n$ are linearly independent 
vectors of $\Sigma_2$ and 
$$\Vert {\bf a}_1,\cdots ,{\bf a}_n\Vert =\sqrt{
\langle {\bf a}_1,\cdots ,{\bf a}_n\vert {\bf a}_1,\cdots ,{\bf a}_n\rangle},
\quad 
\Vert {\bf b}_1,\cdots ,{\bf b}_n\Vert =\sqrt{
\langle {\bf b}_1,\cdots ,{\bf b}_n\vert {\bf b}_1,\cdots ,{\bf b}_n\rangle}.
$$
The angle $\varphi$ does not depend on the choice of the 
bases ${\bf a}_1,\cdots ,{\bf a}_n$ and ${\bf b}_1,\cdots ,{\bf b}_n$. 

Note that any $n$-inner product induces an ordinary inner product 
over the vector space $\Lambda _n(V)$ of $n$-forms on $V$ as follows. 
Let $\{ {\bf e}_{\alpha}\}$, be a basis of $V$. Then we define 
$$\langle \sum_{i_1,\cdots ,i_n}a_{i_1\cdots i_n}
{\bf e}_{i_1}\wedge \cdots \wedge {\bf e}_{i_n}\vert 
\sum_{j_1,\cdots ,j_n}b_{j_1\cdots j_n}
{\bf e}_{j_1}\wedge \cdots \wedge {\bf e}_{j_n}\rangle = $$
$$=\sum_{i_1,\cdots ,i_n,j_1,\cdots ,j_n}
a_{i_1\cdots i_n}b_{j_1\cdots j_n}
\langle {\bf e}_{i_1},\cdots ,{\bf e}_{i_n}\vert 
{\bf e}_{j_1},\cdots ,{\bf e}_{j_n}\rangle .$$
The first requirement for the inner product is a consequence of the 
Theorem 2.1. For example, if 
$$w=p{\bf e}_{i_1}\wedge \cdots \wedge {\bf e}_{i_n}-
q{\bf e}_{j_1}\wedge \cdots \wedge {\bf e}_{j_n},$$ 
then 
$$\langle w\vert w\rangle = 
p^2 \langle {\bf e}_{i_1},\cdots ,{\bf e}_{i_n}\vert 
{\bf e}_{i_1},\cdots ,{\bf e}_{i_n}\rangle + 
q^2 \langle {\bf e}_{j_1},\cdots ,{\bf e}_{j_n}\vert 
{\bf e}_{j_1},\cdots ,{\bf e}_{j_n}\rangle - 
2pq\langle {\bf e}_{i_1},\cdots ,{\bf e}_{i_n}\vert 
{\bf e}_{j_1},\cdots ,{\bf e}_{j_n}\rangle \ge 0$$
and moreover, the last expression is 0 if and only if 
$$\langle {\bf e}_{i_1},\cdots ,{\bf e}_{i_n}\vert 
{\bf e}_{j_1},\cdots ,{\bf e}_{j_n}\rangle =
\sqrt{\langle {\bf e}_{i_1},\cdots ,{\bf e}_{i_n}\vert 
{\bf e}_{i_1},\cdots ,{\bf e}_{i_n}\rangle }
\sqrt{\langle {\bf e}_{j_1},\cdots ,{\bf e}_{j_n}\vert 
{\bf e}_{j_1},\cdots ,{\bf e}_{j_n}\rangle }$$
which means that ${\bf e}_{i_1},\cdots ,{\bf e}_{i_n}$ and 
${\bf e}_{j_1},\cdots ,{\bf e}_{j_n}$ generate the same subspace, and 
$p{\bf e}_{i_1}\wedge \cdots \wedge {\bf e}_{i_n}=
q{\bf e}_{j_1}\wedge \cdots \wedge {\bf e}_{j_n}$, i.e. if and only if 
$w=0$. The other requirements for inner product are obviously satisfied. 
Hence we obtain an induced ordinary inner product on the vector 
space $\Lambda _n(V)$ of $n$-forms on $V$. 

\smallskip 
\noindent 
{\bf Remark.} Note that the inner product on $\Lambda _n(V)$ introduced 
in Example 2.1 is only a special case of inner product on 
$\Lambda _n(V)$ and also $n$-inner product. It is induced via 
existence of an ordinary inner product on $V$. 

The angle between 
subspaces defined by (3.1) coincides with the angle between two $n$-forms 
in the vector space $\Lambda _n(V)$. Since the angle between two "lines"
in any vector space with ordinary inner product can 
be considered as a distance, we obtain that 
$$\varphi = \arccos {
\langle {\bf a}_1,\cdots ,{\bf a}_n\vert {\bf b}_1,\cdots ,{\bf b}_n\rangle 
\over 
\Vert {\bf a}_1,\cdots ,{\bf a}_n\Vert \cdot 
\Vert {\bf b}_1,\cdots ,{\bf b}_n\Vert }\leqno{(3.2)}$$
determines a metric among the $n$-dimensional subspaces of $V$. Indeed, 
it induces a metric on the Grassmann manifold $G_n(V)$, which is 
compatible with the ordinary topology of the Grassman manifold $G_n(V)$. 
This metric over Grassmann manifolds appears natural 
and appears convenient also for the infinite dimensional vector spaces $V$. 

Further we shall consider a special case of $n$-inner product for which 
there exists a basis $\{ {\bf e}_{\alpha}\}$ of $V$ such that the vector 
${\bf e}_i$ is orthogonal to the subspace generated by the vectors 
${\bf e}_{i_1},\cdots ,{\bf e}_{i_n}$ for different values of 
$i,i_1,\cdots ,i_n$. For such $n$-inner product we have 
$$\langle {\bf e}_{i_1},\cdots ,{\bf e}_{i_n}\vert {\bf e}_{j_1}\cdots 
{\bf e}_{j_n}\rangle = 
C_{i_1\cdots i_n}\delta ^{i_1\cdots i_n}_{j_1\cdots j_n}\leqno{(3.3)}$$
where 
$\delta ^{i_{1}\cdots i_{n}}_{j_{1}\cdots j_{n}}$ is equal to 1 or -1 
if $\{ i_1,\cdots ,i_n\} = \{ j_1,\cdots ,j_n\}$ with different 
$i_1,\cdots ,i_n$ and additionally the permutation 
${i_1\; i_2\; \cdots \;i_n\choose j_1\; j_2\; \cdots \;j_n}$ 
is even or odd respectively, and the expression is 0 otherwise, 
and where $C_{i_1\cdots i_n}>0$. Moreover, one can verify that  the  previous 
formula induces an $n$-inner product, i.e. the six conditions i) - vi) are 
satisfied if and only if all the coefficients $C_{i_1\cdots i_n}$ are equal 
to a  positive constant $C>0$. Moreover, we can assume that 
$C=1$, because otherwise we can consider the basis 
$\{ {\bf e}_{\alpha}/C^{1/2n}\}$ instead of the basis 
$\{ {\bf e}_{\alpha}\}$ of $V$. Hence this special case of $n$-inner product
reduces to the $n$-inner product given by the Example 2.1. Indeed, 
the ordinary inner product is uniquely defined such that 
$\{ {\bf e}_{\alpha}\}$ as an orthonormal system of vectors. 

If the dimension of $V$ is finite, for example $\dim V=m>n$, 
then the previous $n$-inner product induces a dual $(m-n)$-inner product on 
$V$ which is induced by 
$$\langle {\bf e}_{i_1},\cdots ,{\bf e}_{i_{m-n}}\vert {\bf e}_{j_1}\cdots 
{\bf e}_{j_{m-n}}\rangle^* = 
\delta^{i_1\cdots i_{m-n}}_{j_1\cdots j_{m-n}}.\leqno{(3.4)}$$
The dual $(m-n)$-inner product is defined using the "orthonormal basis"
$\{ {\bf e}_{\alpha}\}$ of $V$. If we have chosen another "orthonormal 
basis", the result will be the same. Further we prove the following theorem. 

\smallskip 
\noindent 
{\bf Theorem 3.1.} {\em Let $V$ is a finite dimensional vector space 
and let the $n$-inner product on $V$ be defined as in the Example 2.1. Then
$$ \varphi (\Sigma _{1},\Sigma _{2}) = 
\varphi (\Sigma _{1}^{*},\Sigma _{2}^{*}),$$
where $\Sigma_1$ and $\Sigma_2$ are arbitrary $n$-dimensional subspaces of 
$V$ and $\Sigma_1^*$ and $\Sigma_2^*$ are their orthogonal subspaces in $V$.}

\smallskip
\noindent {\em Proof.} Let $\Sigma _{1}=\langle \omega _{1}\rangle$, 
$\Sigma _{2}=\langle \omega _{2}\rangle$, 
$\Sigma _{1}^{*}=\langle \omega _{1}^{*}\rangle$, $\Sigma _{2}^{*}=\langle 
\omega _{2}^{*}\rangle$, 
where $\Vert \omega _{1}\Vert =\Vert \omega _{2}\Vert = 
\Vert \omega _{1}^{*}\Vert =\Vert \omega _{2}^{*}\Vert =1$. We will prove 
that 
$$ \omega _{1}\cdot \omega _{2} = \pm \omega _{1}^{*}\cdot \omega _{2}^{*}.$$
Indeed, 
$\omega _{1}\cdot \omega _{2} = \omega _{1}^{*}\cdot \omega _{2}^{*}$ if 
$\omega _{1}\wedge \omega _{2}$ and $\omega _{1}^{*}\wedge \omega _{2}^{*}$
have the same orientation in $V$ and 
$\omega _{1}\cdot \omega _{2} = -\omega _{1}^{*}\cdot \omega _{2}^{*}$ if 
$\omega _{1}\wedge \omega _{2}$ and $\omega _{1}^{*}\wedge \omega _{2}^{*}$
have the opposite orientations in $V$. 
\par 
Assume that the dimension of $V$ is $m$. 
Without loss of generality we can assume that 
$$ \omega _{1}={\bf e}_{1}\wedge {\bf e}_{2}\wedge \cdots \wedge {\bf e}_{n}
\quad \hbox { and }\quad 
\omega _{1}^{*}=
{\bf e}_{n+1}\wedge {\bf e}_{n+2}\wedge \cdots \wedge {\bf e}_{m}.$$
Without loss of generality we can assume that 
$$ \omega _{2}={\bf a}_{1}\wedge {\bf a}_{2}\wedge \cdots \wedge {\bf a}_{n}
\quad \hbox { and }\quad 
\omega _{2}^{*}=
{\bf a}_{n+1}\wedge {\bf a}_{n+2}\wedge \cdots \wedge {\bf a}_{m},$$
where ${\bf a}_{1}, \cdots ,{\bf a}_{m}$ is an orthonormal system. 
Suppose that 
${\bf a}_{i}=(a_{i1},\cdots ,a_{im})$ $(1\le i\le m)$, and let us introduce 
an orthogonal $m\times m$ matrix 
$$ 
A = \left [\matrix{ 
a_{11}& \cdots &a_{1n} & a_{1,n+1}& \cdots &a_{1m}\cr
\cdot &  &  &  &  &  \cr
\cdot &  &  &  &  &  \cr
\cdot &  &  &  &  &  \cr
a_{n1}& \cdots &a_{nn} & a_{n,n+1}& \cdots &a_{nm}\cr
a_{n+1,1}& \cdots &a_{n+1,n} & a_{n+1,n+1}& \cdots &a_{n+1,m}\cr
\cdot &  &  &  &  &  \cr
\cdot &  &  &  &  &  \cr
\cdot &  &  &  &  &  \cr
a_{m1}& \cdots &a_{mn} & a_{m,n+1}& \cdots &a_{mm}\cr }\right ] .$$
We denote by $A_{i_{1}\cdots i_{n}}$ $(1\le i_{1}<i_{2}<\cdots <i_{n}\le m)$,
the $n\times n$ submatrix of $A$ whose rows are the first $n$ rows of $A$
and whose columns are the $i_{1}$-th,...,$i_{n}$-th column of $A$. 
We denote by $A^{*}_{i_{1}\cdots i_{n}}$ the $(m-n)\times (m-n)$ submatrix 
of $A$ which is obtained 
by deleting the rows and the columns corresponding to 
the submatrix $A_{i_{1}\cdots i_{n}}$. It is easy to verify that 
$$ \omega _{1}\cdot \omega _{2} = det A_{12 ... n} \quad \hbox {and }\quad 
\omega _{1}^{*}\cdot \omega _{2}^{*} = det A_{12 ... n}^{*}$$
and thus we have to prove that 
$$ det A_{12 ... n} = \pm det A_{12 ... n}^{*}, \leqno{(3.5)}$$
i.e. 
$$ det A_{12 ... n} = det A_{12 ... n}^{*}\quad \hbox { if }\quad det A=1$$
and 
$$det A_{12 ... n} = -det A_{12 ... n}^{*}\quad \hbox { if }\quad det A=-1.$$

Assume that $det A=1$. Let us consider the expression 
$$ F = 
\sum _{1\le i_{1}<i_{2}<\cdots <i_{n}\le m}\Bigl [
(det A_{i_{1}i_{2}\cdots i_{n}} - (-1)^{1+2+\cdots +n}
(-1)^{i_{1}+i_{2}+\cdots +i_{n}}
det A_{i_{1}i_{2}\cdots i_{n}}^{*})\Bigr ]^{2}.$$
Using that $\Vert \omega _{2}\Vert=1$ and $\Vert \omega _{2}^{*}\Vert=1$
we get 
$$\sum _{1\le i_{1}<\cdots <i_{n}\le m} 
(det A_{i_{1}i_{2}\cdots i_{n}})^{2} = 
\sum _{1\le i_{1}<\cdots <i_{n}\le m} 
(det A_{i_{1}i_{2}\cdots i_{n}}^{*})^{2} = 1$$
and using the Laplace formula for decomposition of determinants, we obtain 
$$ F = 
\sum _{1\le i_{1}<\cdots <i_{n}\le m} 
(det A_{i_{1}i_{2}\cdots i_{n}})^{2} + 
\sum _{1\le i_{1}<\cdots <i_{n}\le m} 
(det A_{i_{1}i_{2}\cdots i_{n}}^{*})^{2}  $$
$$ -2 
\sum _{1\le i_{1}<\cdots <i_{n}\le m} (-1)^{n(n+1)/2}
(-1)^{i_{1}+i_{2}+\cdots +i_{n}}
det A_{i_{1}i_{2}\cdots i_{n}}
det A_{i_{1}i_{2}\cdots i_{n}}^{*} $$
$$ = 1+1-2\cdot det A = 2 - 2 = 0 .$$
Hence $F=0$ implies that 
$$ det A_{i_{1}i_{2}\cdots i_{n}} = (-1)^{n(n+1)/2}
(-1)^{i_{1}+i_{2}+\cdots +i_{n}}det A_{i_{1}i_{2}\cdots i_{n}}^{*}.$$
Specially for $i_{1}=1, \cdots ,i_{n}=n$ we obtain 
$$ det A_{i_{1}i_{2}\cdots i_{n}} = det A_{i_{1}i_{2}\cdots i_{n}}^{*}.$$
\par 
Assume that $det A=-1$. Then we consider the expression 
$$ F' = 
\sum _{1\le i_{1}<i_{2}<\cdots <i_{n}\le m}\Bigl [
(det A_{i_{1}i_{2}\cdots i_{n}} + (-1)^{1+2+\cdots +n}
(-1)^{i_{1}+i_{2}+\cdots +i_{n}}
det A_{i_{1}i_{2}\cdots i_{n}}^{*})\Bigr ]^{2}$$
and analogously we obtain that 
$$ det A_{i_{1}i_{2}\cdots i_{n}} = -(-1)^{n(n+1)/2}
(-1)^{i_{1}+i_{2}+\cdots +i_{n}}det A_{i_{1}i_{2}\cdots i_{n}}^{*}.$$
Specially for $i_{1}=1, \cdots ,i_{n}=n$ we obtain 
$$ det A_{i_{1}i_{2}\cdots i_{n}} = -det A_{i_{1}i_{2}\cdots i_{n}}^{*}.
\eqno{\Box }$$

Finally we remark the following. The presented approach to the $n$-inner 
product appears essential for the applications in the functional analysis. 
Since the corresponding $n$-norm is the same as the corresponding $n$-norm 
from the definition of Misiak, we have the same results in the normed 
spaces. It is an open question whether from the Definition 2.1 
can be introduced a generalized $n$-inner product and $n$-semiinner product 
with characteristic $p$. It is of interesting also to research the strong 
convexity in the possibly introduced space with $n$-semiinner product
with characteristic $p$.

\noindent ${}^1$Institute of Mathematics

\noindent Ss. Cyril and Methodius University in Skopje

\noindent P.O.Box 162, 1000 Skopje, Macedonia

\noindent e-mail: kostatre@pmf.ukim.mk, kostadin.trencevski@gmail.com
\medskip

\noindent ${}^2$Faculty of Social Sciences 

\noindent Anton Popov, b.b. 

\noindent 1000 Skopje, Macedonia

\noindent e-mail: risto.malceski@gmail.com  

\begin{thebibliography}{9}
\bibitem{F} S. FEDOROV: Angle between subspaces of analytic and 
antianalytic functions in weighted $L_2$ space on a boundary of a multiply
connected domain, in: {\em Operator Theory. System Theory and Related 
Topics}, Beer-Sheva/Rehovot (1997), 229-256. 

\bibitem{Kn} A. V. KNYAZEV and M. E. ARGENTATI: Principal angles between 
subspaces in an A-based scalar product: algorithms and perturbation 
estimates, {\em SIAM J. Sci. Comput.}, {\bf 23} (2002), 2008-2040. 

\bibitem{K} S. KUREPA: On the Buniakowsky-Cauchy-Schwarz inequality, 
{\em Glasnik Mat. Ser. III}, {\bf 1}(21) (1966), 147-158. 

\bibitem{1} {\sc S. MACLANE and G. BIRKHOFF}: Algebra, 
{\em The Macmillan Company}, New York (1967). 

\bibitem{2} {\sc A. MISIAK}: 
{\em n}-Inner Product Spaces, 
{\em Math. Nachr.}, {\bf 140} (1989), 299-319. 

\bibitem{RW} V. RAKO\v{C}EVI\'{C} and H. K. WIMMER: 
A variational characterization of canonical angles between subspaces, 
{\em J. Geom.}, {\bf 78} (2003), 122-124. 

\bibitem{W} H. K. WIMMER: Canonical angles of unitary spaces and 
perturbations of direct complements, {\em Linear Algebra Appl.}, {\bf 287} 
(1999), 373-379. 


\end{thebibliography}
\end{document}